\documentclass[12pt]{amsart}
\usepackage{amssymb,latexsym,amsmath,amsthm,enumitem,mathrsfs,geometry}
\usepackage{fullpage}
\usepackage{fancyhdr}
\usepackage{xcolor}
\usepackage{bbm}
\usepackage{stmaryrd}
\usepackage{mathrsfs}
\usepackage{hyperref}
\usepackage{lineno}
\usepackage{diagbox}
\usepackage{array}
\usepackage{tikz}
\usetikzlibrary{arrows}
\usetikzlibrary{positioning, arrows.meta, shapes.geometric}
\usepackage{graphicx}
\usepackage{float}
\usepackage{comment}
\usepackage{mathtools}
\usepackage{cleveref}
\usepackage{tikz-3dplot}

\newtheorem*{theorem*}{Theorem}
\newtheorem{problem}{Problem}[section]
\newtheorem{conjecture}{Conjecture}[section]

\newtheorem*{remark*}{Remark}

\newcommand{\eps}{\varepsilon}
\newcommand{\R}{\mathbb{R}}
\newcommand{\C}{\mathbb{C}}

\newcommand{\Z}{\mathbb{Z}}

\begin{document}

\numberwithin{equation}{section}

\title{Rotating needles in space: the road to the Kakeya conjecture, and why it matters}

\author{Terence Tao}
\address{UCLA Department of Mathematics, Los Angeles, CA 90095-1555.}
\email{tao@math.ucla.edu}

\thanks{}

\date{\today}

\begin{abstract} A non-technical exposition of the Kakeya conjecture, why it matters, and the road to the solution of this conjecture in three dimensions by Hong Wang and Joshua Zahl.
\end{abstract}

\maketitle

\section{The Kakeya problem}

As any car driver knows, it is possible to maneuver a car into a tight parking spot, or to perform a U-turn in a narrow street, by performing a sequence of forward and backward motions.  A mathematician might then be inclined to ask just how much space is needed to execute, say, a U-turn?
This of course depends on how large and wide the car is, so the mathematician might first work in an idealized setting in which the car has zero thickness (i.e., is a line segment or ``needle''), and is infinitely maneuverable.  This leads to the following question, which was first posed by Soichi Kakeya \cite{kakeya} in 1917:

\begin{problem}[Kakeya needle problem]\label{needle} What is the smallest area of a planar region in which a unit line segment can be continuously rotated by a full rotation\footnote{Strictly speaking, a $U$-turn would only consist of a half-rotation by $180^\circ$, but the two questions are very similar, and the needle problem is traditionally formulated using a full rotation.} of $360^\circ$?
\end{problem}

\begin{figure}[h]
\centering
\includegraphics[width=0.6\textwidth]{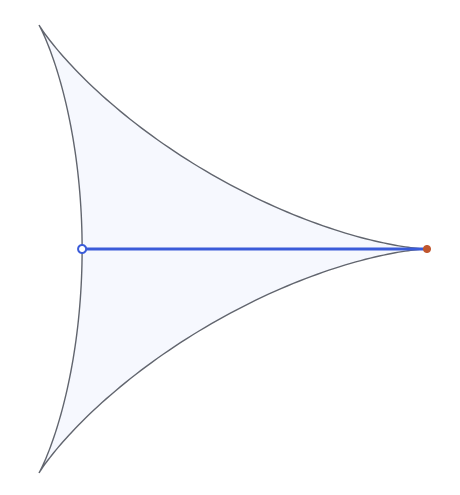}
\caption{A deltoid of area $\pi/8$ in which a unit needle (the blue line segment) can be rotated.}\label{deltoid}
\end{figure}

For instance, if one turns the needle around its midpoint, it sweeps out a disk of radius $1/2$, which has area $\pi/4 \approx 0.785$.  A ``three-point $U$-turn'' can replace this disk by a smaller shape, known as a \emph{deltoid}, which has area $\pi/8 \approx 0.393$; see Figure \ref{deltoid}.

\begin{figure}[h]
\centering
\includegraphics[width=0.6\textwidth]{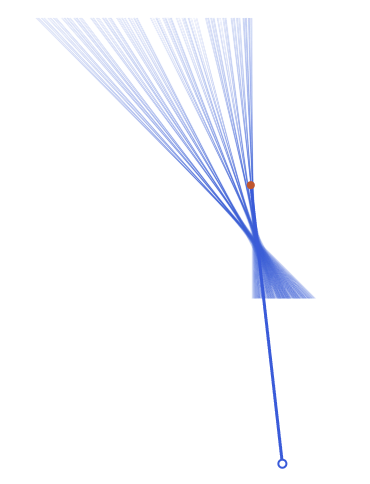}
\caption{A Besicovitch--Perron tree, which contains a unit line segment in every direction in an $45^\circ$ arc, but can be chosen to have arbitrary small area.  It is possible to make a needle perform a continuous $U$ turn using four such trees, together with some very small additional strips (not pictured) to enable the needle to move from one branch of the tree to the next. An animated applet demonstrating this tree may be found at \url{https://teorth.github.io/tao-web/apps/kakeya.html}.}\label{tree}
\end{figure}

For a while, it was believed that this was the most efficient solution to the needle problem; but in 1928, Abram Besicovitch \cite{besicovitch} made the unexpected discovery that one can in fact rotate such a needle using arbitrarily small area.  A key ingredient was the construction of what is now known as a \emph{Besicovitch--Perron tree}, an example of which is given in Figure \ref{tree}.

The Besicovitch--Perron tree is (one quarter of) a type of fractal set now known as a \emph{Besicovitch set} or a \emph{Kakeya set} - a set that contains a unit line segment in every direction.  Besicovitch showed that such sets can be constructed to have an area as small as one pleases; in fact, one can even find Besicovitch sets of area zero\footnote{There is a technical subtlety: despite the existence of Besicovitch sets of area zero, for the original needle problem one still needs to add a small additional amount of area to allow the needle to be rotated continuously.  So the full solution to Problem \ref{needle} is ``as small as one pleases, as long as it is non-zero''.}!  We will not give the details of these constructions here, save to say that like many other fractal sets, such Besicovitch sets have a multi-scale structure, where pairs of line segments making small, medium, or large angles to each other are all pushed together strategically to create a significant amount of overlap in area.

Now suppose that the car is not an infinitely thin needle, but instead has some positive thickness.  Mathematically, one can model this by a rectangle of dimensions $1 \times \delta$ of some small $\delta > 0$.  Now of course one cannot possibly hope to rotate such a car using arbitrarily small area, as one must, at a minimum, have enough area to accommodate the car itself, which now has area $\delta$.  If one takes a Besicovitch set and forms its \emph{$\delta$-neighborhood} - the set of all points within a distance $\delta$ of a that set - then that creates a set large enough to accommodate every rotation of the $1 \times \delta$ rectangle.  The natural question is then: how large must that $\delta$-neighborhood be? Or, to put it another way, if one has a copy of a $1 \times \delta$ in every possible orientation, and one is allowed to slide (but not rotate) these rectangles around to make them overlap, what is the smallest area one can compress these rectangles into?

This turns out to be a more delicate question, with satisfactory answers only coming with the work of C\'ordoba in 1977 \cite{cordoba} and Keich in 1999 \cite{keich}.  Combining their results, they showed that the area of such a neighborhood shrinks to zero as $\delta$ goes to zero, but only at the very slow logarithmic rate of about\footnote{We will be vague here about what ``about'' means, but it roughly means ``up to multiplicative constants not depending on the thickness $\delta$''.} $1/\log(1/\delta)$.   For instance, a needle of thickness $\delta = 10^{-100}$ would still need a region of area about $1/100$ to rotate in: small, but not nearly as tiny as the area of the needle itself.
Thus, while a certain amount of ``Kakeya compression'' is possible, the extent of such compression is sharply limited.

Why is it so hard to compress?  A back of the envelope calculation is as follows.  One needs to rotate a $1 \times \delta$ rangle by an angle of at least $\delta$ to get a genuinely new orientation of the rectangle (one that can't be easily covered by one or two translates of the original rectangle).  So there are roughly about $1/\delta$ different orientations of the rectangle that one needs to compress together.  Typically, the intersection between two such rectangles becomes a parallelogram of dimensions roughly of order $\delta$, leading to an area of about $\delta^2$, or a proportion of about $\delta$ of each of the rectangles.  As there are only $1/\delta$ or so orientations, this suggests that a typical point in one rectangle is covered by only one or two other rectangles, which is not enough overlap to create genuine compression.  This argument is oversimplified because it does not take into account the fact that rectangles that are close in orientation can be made to overlap in an area significantly larger than $\delta^2$; but the aforementioned paper of C\'ordoba \cite{cordoba} accounted for this fact more carefully, and demonstrated that the total compression factor could only grow logarithmically in $\delta$ at best.

\begin{figure}[h]
\centering
\includegraphics[width=0.6\textwidth]{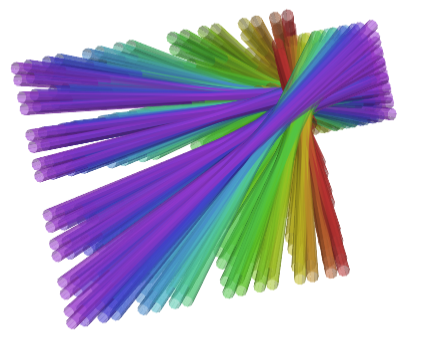}
\caption{A configuration of $1 \times \delta \times \delta$ cylinders pointing in different directions.  An animated visualization of this three-dimensional analogue of a Besicovitch--Perron tree may be found at \url{https://teorth.github.io/tao-web/apps/kakeya3d.html}.}\label{tree3d}
\end{figure}

These results (as well as some more technical variants which I will not discuss here) give a pretty satisfactory answer to Kakeya type problems in two dimensions.  But a mathematician always keeps an eye on possible extensions and generalizations.  One natural question is: what happens in three dimensions?  Now, instead of a thin $1 \times \delta$ rectangle, imagine a thin cylinder (or ``tube'') of dimensions $1 \times \delta \times \delta$ in space, and suppose one wants to find a ``$\delta$-thickened Kakeya set'' that contains a copy of every possible orientation of this tube; see Figure \ref{tree3d}.  (One could perhaps imagine this tube as representing a space telescope that is intended to look at every star in the sky, but for some inexplicable reason needs to be confined to as small a volume of space as possible.)  The question is then: how small can the volume of such a set be?

There are many forms of this problem, some of which are technical to state, but to focus the presentation, we shall just discuss one its forms, introduced by Bourgain \cite{bourgain}:

\begin{conjecture}[Kakeya conjecture]\label{kakeya-conj}  Let $0 < \delta < 1$ be a small quantity, and $E$ be a subset of $\R^3$ that contains a $1 \times \delta \times \delta$ cylinder in every direction.  Must the volume of $E$ be at least $c_\eps \delta^\eps$ for every fixed $\eps>0$?  Here $c_\eps>0$ is a constant that can depend on $\eps$ but not on $\delta$.
\end{conjecture}

Informally, this conjecture asserts that -- as in the two-dimensional case -- the amount of ``compression'' that is possible amongst tubes in three dimensions pointing in different directions can only be (roughly) logarithmic in the eccentricity $1/\delta$ of the cylinder.  Using asymptotic notation, the Kakeya conjecture can be stated more succicntly as the assertion that the volume of $E$ is at least $\delta^{o(1)}$ as $\delta \to 0$.  

\section{Why is the Kakeya problem important?}

The Kakeya problem sounds like a pure curiosity; after all, there isn't much direct application for compressing a space telescope into a small volume.  However, it turns out that the problem of understanding the compressibility of tubes (or tube-like objects) is central to many problems in harmonic analysis, partial differential equations, and number theory, which all have the general flavor of trying to control the amount of ``constructive interference'' that can be present in some superposition of oscillating waves in ``space''.

An early, and particularly striking, example of this connection was the work of Fefferman in 1971 \cite{fefferman} in Fourier analysis, who famously gave a negative answer to the \emph{disk multiplier problem}, which we will now pause to explain.  The fundamental work of Fourier showed, in principle at least, that every (reasonably well behaved) function $f : \R \to \C$ that is $1$-periodic can be expressed as a Fourier series
$$
f(x) = \sum_{n \in \Z} \hat f(n) e^{2\pi i n x}$$
for some coefficients $\hat f(n) \in \C$.  For instance, the square wave function
$$ f(x) = \begin{cases} 1 & \text{if } 0 \leq \{x\} < 1/2 \\ -1 & \text{if } 1/2 \leq \{x\} < 1 \end{cases},$$
where $\{x\}$ denotes the fractional part of $x$, has a Fourier series expansion
$$
f(x)=\sum_{n \text{ odd}} \frac{-2i}{\pi n} e^{2\pi i n x}.
$$
Joseph Fourier introduced his famous series expansions in 1807 in order to solve the heat equation.  Fourier's series expansions were initially controversial, as his calculations implied that piecewise discontinuous functions, such as the square wave, could be expressed as an infinite series of continuous sinusoids, which went against the mathematicial intuition of the time.  The matter was not fully settled until much later, when the methods of analysis had been put on a much firmer footing.  In modern language, the basic question is to determine the senses in which the
\emph{partial Fourier series}
$$
S_N f(x) \coloneqq \sum_{n=-N}^N \hat f(n) e^{2\pi i n x}$$
converge back to the original function.

\begin{figure}[h]
\centering
\includegraphics[width=0.6\textwidth]{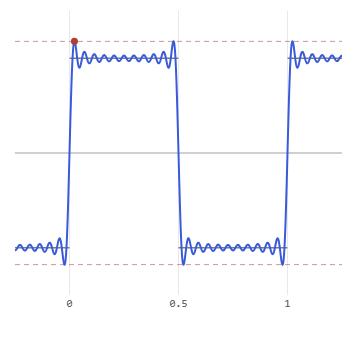}
\caption{The Gibbs phenomenon for the square wave function.}\label{gibbs}
\end{figure}

The answer to this question is surprisingly subtle, and depends on how smooth $f$ is.  For functions $f$ that are smooth enough (continuously differentiable will do), the partial Fourier series converge uniformly to $f(x)$; but in the presence of jump discontinuities one has the infamous \emph{Gibbs phenomenon}, in which the Fourier series stay a fixed distance from the function near the jump, no matter how large one takes $N$; see Figure \ref{gibbs}.

Nevertheless, it follows from the foundational work of F. Riesz and Fischer in 1907 \cite{riesz-1907, fischer} (or the later theorem of Plancherel \cite{plancherel} in 1910) that if $f$ is in the square-integrable class $L^2$, then the Fourier series also converge in $L^2$, thus
$$ \int_0^1 |S_N f(x) - f(x)|^2\ dx \to 0$$
as $N \to \infty$.  Informally, this is an assertion that the type of pathologies seen in the Gibbs phenomenon are localized in space, and become asymptotically negligible upon performing a square-averaging in space.  This result was later extended to the $L^p$ setting by M. Riesz in 1928 \cite{riesz}: if $\int_0^1 |f(x)|^p\ dx < \infty$ for some\footnote{The endpoint cases $p=1$ and $p=\infty$ of this claim are false, for reasons we will not go into here.} $1 < p < \infty$, then
$$ \int_0^1 |S_N f(x) - f(x)|^p\ dx \to 0$$
as $N \to \infty$.  

One can think of the partial Fourier series $S_N f$ as a \emph{Fourier multiplier} of $f$, because the Fourier coefficients of $S_N f$ are obtained from those of $f$ by multiplication:
$$ \widehat{S_N f}(n) = 1_{[-N,N]}(n) \hat f(n).$$

In the language of harmonic analysis, the convergence theorems of F. Riesz, Fischer, and M. Riesz are basically asserting that the indicator functions $1_{[-N,N]}$ of intervals $[-N,N]$ are \emph{$L^p$ multipliers} uniformly in $N$, for every fixed $1 < p < \infty$.  The theory of $L^p$ multipliers, both in one dimension and in higher dimensions, flourished for much of the twentieth century, with many fundamental \emph{multiplier theorems} by Marcinkiewicz, H\"ormander--Mikhlin, and others forming the backbone of modern harmonic analysis and the analysis of partial differential equations, leading for instance to the foundational theory of Sobolev spaces and pseudo-differential operators.  Roughly speaking, these theorems allowed one to go beyond the $L^2$-based Fourier analysis provided by the Plancherel theorem, and allow one to manipulate Fourier series or Fourier integrals in other $L^p$ spaces.

However, once one went to higher dimensions, there were some basic Fourier multipliers that strangely resisted all attempts to tame them.  The simplest example would be the \emph{disk multipliers} $D_N$ in two dimensions.  A function $f(x,y)$ of two variables that was periodic with period $1$ in each variable would have a two-dimensional Fourier series
$$
f(x,y) = \sum_{(m,n) \in \Z^2} \hat f(m,n) e^{2\pi i (m x + n y)}.$$
If one decided to sum this series over squares of frequencies $\{ (m,n): -N \leq m,n \leq N \}$ to create square multipliers
$$ S_N f(x,y) \coloneqq \sum_{-N \leq m,n \leq N} \hat f(m,n) e^{2\pi i (m x + n y)},$$
then one could use the one-dimensional multiplier theory to show that these partial Fourier series converged back to the original function $f$ as $N \to \infty$ in the $L^p$ sense for any $1 < p < \infty$, assuming of course that the original function $f$ was in the $L^p$ class.  But if instead one made the innocuous change of replacing the square $\{ (m,n): -N \leq m,n \leq N \}$ by the disk $\{ (m,n): m^2+n^2 \leq N^2 \}$, then the resulting disk multipliers
$$ D_N f(x,y) \coloneqq \sum_{m^2+n^2 \leq N^2} \hat f(m,n) e^{2\pi i (m x + n y)}$$
did not seem to be controllable by any existing theory, (except in the $p=2$ case, where the theorem of Plancherel easily gave convergence).  

\begin{figure}[h]
\centering
\includegraphics[width=0.6\textwidth]{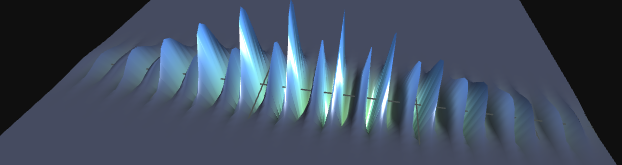}
\caption{A gaussian wave packet, oscillating at a certain frequency while also being localized in space.}\label{wavepacket}
\end{figure}

It came as a great surprise when Fefferman in 1977 showed that in fact the disk multipliers $D_N f$ did \emph{not} always converge in the $L^p$ sense back to the original function $f$ for $p \neq 2$.  The result was shocking enough, but what was more surprising at the time was the \emph{method} used: not the functional analysis or operator algebra techniques that were ascendant at the time, but instead a geometric construction closely related to the Kakeya problem.  In modern language, the idea was this.  Firstly, through some basic Fourier analytic computations, one could generate a \emph{wave packet}: an oscillating function concentrated on a thin rectangle in the plane (and then extended periodically), whose Fourier coefficients were just on the edge of the disk $\{ (m,n): m^2+n^2 \leq N^2 \}$; see Figure \ref{wavepacket}.  If the parameters were chosen appropriately, applying the disk multiplier would only retain about half of the frequencies comprising the wave packet, which (by the Heisenberg uncertainty principle) would cause the wave packet to elongate, filling out a slightly larger rectangle in the plane than it initially lived in.  By modifying the Bescovitch--Perron construction, Fefferman was able to arrange a number of such wave packets in such a way that initially they were separated from each other and did not exhibit much interference, but upon lengthening the wave packets, they enter a compressed configuration, reminiscent of a Kakeya set, in which the wave packets interfere with each other and create a large amplitude.  This formed the basis of Fefferman's counterexample of a Fourier series in which the disk multipliers failed to converge in the expected sense.  For more details, see the recent survey of Hickman \cite{hickman}.

In the decades that followed, many further connections of this type were unearthed, in which the objects being analyzed contained ``wave packets'' that could point in different directions due to some ``curvature'' in the problem could end up being compressed into a Kakeya-like configuration, thus leading to unwanted behavior.  Because of this, many long-standing conjectures in harmonic analysis, number theory, and partial differential equations were shown to imply the Kakeya conjecture, in that any counterexample to the Kakeya conjecture could be converted to a counterexample to these latter conjectures.
Some examples include:

\begin{itemize}
\item \emph{Stein's restriction conjecture}, which describes when the Fourier transform of a multidimensional function can be restricted to a curved surface such as a sphere, cone, or paraboloid.
\item The \emph{Bochner--Riesz conjecture}, which describes the extent to which the failure of the disk multiplier (or higher dimensional analogues) to have good convergence behavior can be repaired by smoothing out the multiplier near the boundary of the disk.
\item Sogge's \emph{local smoothing conjecture}, which describes the extent to which solutions to the wave equation can be smoothed out by averaging in time.
\item The \emph{Montgomery conjecture}, which describes the extent to which generalizations of the Riemann zeta function can attain unusually large values.
\end{itemize}

A more detailed discussion of all of these conjectures, and their relation to the Kakeya problem, may be found in \cite{wolff}.  Here, we shall just explore the Montgomery conjecture, which was an unexpected relation between the Kakeya problem and analytic number theory.  (For more discussion of the other conjectures in this list, see my earlier article \cite{tao-notices}.)

\begin{figure}[h]
\centering
\includegraphics[width=0.9\textwidth]{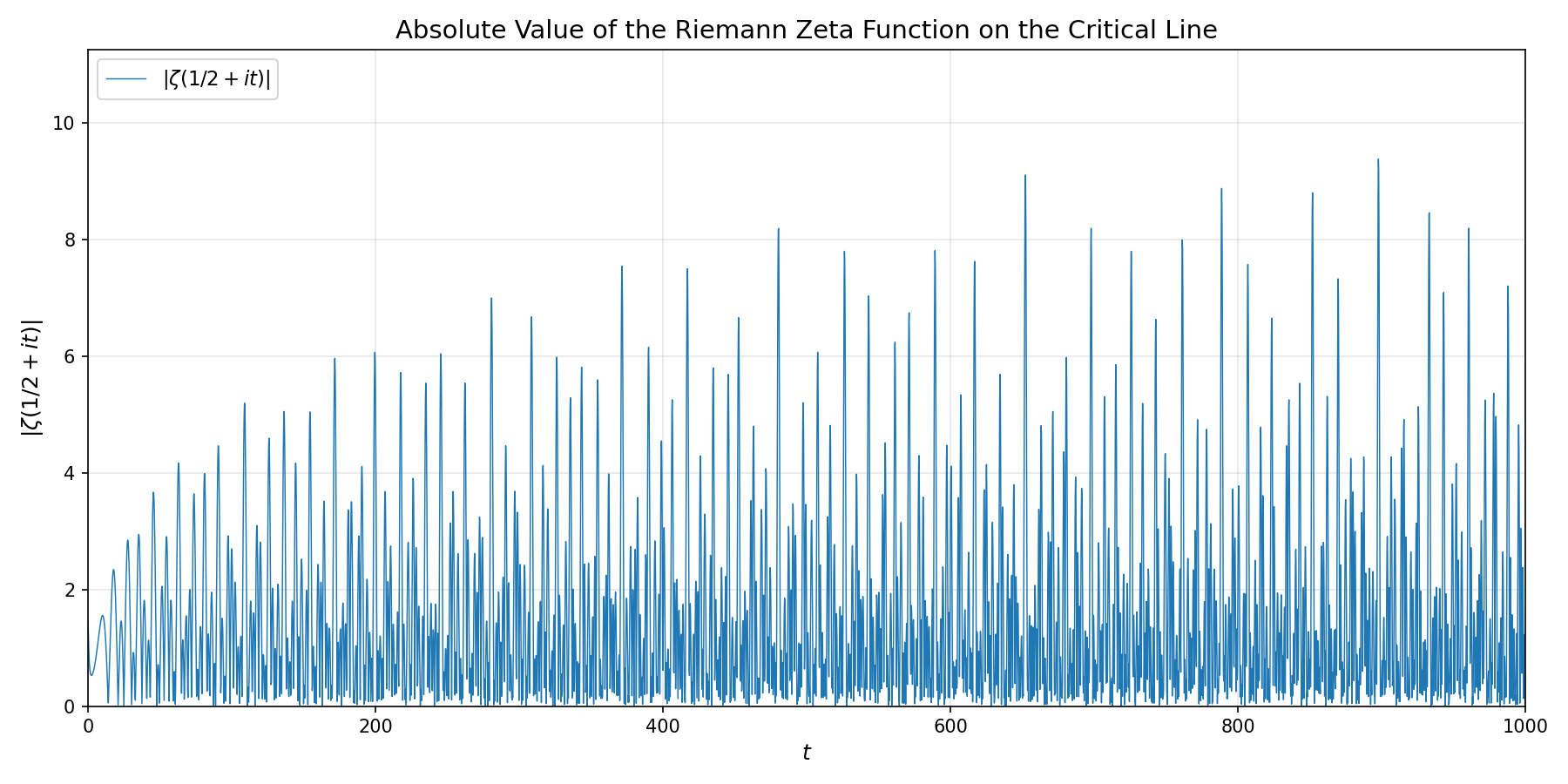}
\caption{The (highly oscillatory) function $|\zeta(\frac{1}{2}+it)|$ is observed numerically to grow quite slowly as $t \to \infty$, which is consistent with the Lindel\"of (or Riemann) hypothesis.}\label{zeta}
\end{figure}

The behaviour of the prime numbers is closely tied to the \emph{Riemann zeta function} $\zeta(s)$, defined (when the real part of $s$ is at least $1$) by the \emph{Dirichlet series}
$$ \zeta(s) \coloneqq \sum_{n=1}^\infty \frac{1}{n^s}.$$
Indeed, in Riemann's famous 1859 memoir \cite{riemann}, the distribution of the prime numbers was shown to be closely related to the distribution of the zeros of $\zeta(s)$.  Through the machinery of complex analysis (such as Jensen's formula), one can relate this distribution to the growth of the zeta function, particularly on the critical line $s = \frac{1}{2} + it$.  For instance, the infamous Riemann hypothesis has long been known to imply the \emph{Lindel\"of hypothesis} that asserts the subpolynomial growth bound $|\zeta(\frac{1}{2}+it)| \leq t^{o(1)}$ as $t \to \infty$; see Figure \ref{zeta}.  This hypothesis, like the Riemann hypothesis, remains open; but the question of growth of this zeta function remains a major focus of study in analytic number theory.

The Lindel\"of hypothesis can also be phrased in terms of the \emph{Dirichlet polynomials}
$$ \sum_{n=1}^N n^{it}$$
where $N$ and $t$ are large parameters; for instance, the hypothesis is equivalent to such polynomials growing no faster than $t^{1/2+o(1)}$ when $N = t^{1/2}$.  Such a uniform bound is not known; however, \emph{mean value theorems} were developed over the years to establish that such bounds are true for \emph{most} values of $t$, rather than \emph{all} values. (This in turn leads to useful partial results towards the Riemann hypothesis, in which most (but not necessarily all) of the zeroes are shown to be close to the critical line.)

These mean value theorems had the feature of generalizing to other Dirichlet polynomials
\begin{equation}\label{dirichlet-polynomial}
\sum_{n=1}^N a_n n^{it}
\end{equation}
where the coefficients $a_n$ were bounded in size, but otherwise arbitrary; thus, they were not specific to the Riemann zeta function, but instead arose from more general inequalities in analysis that had little to do with number theory.  Montgomery conjectured a generalization of many of the existing mean value theorems that would be valid for such arbitrary coefficients, and which would imply the Lindel\"of conjecture, thus opening up the tantalizing possibility of making a major breakthrough towards the Riemann hypothesis primarily by analytic means.

\begin{figure}[h]
\centering
\includegraphics[width=0.9\textwidth]{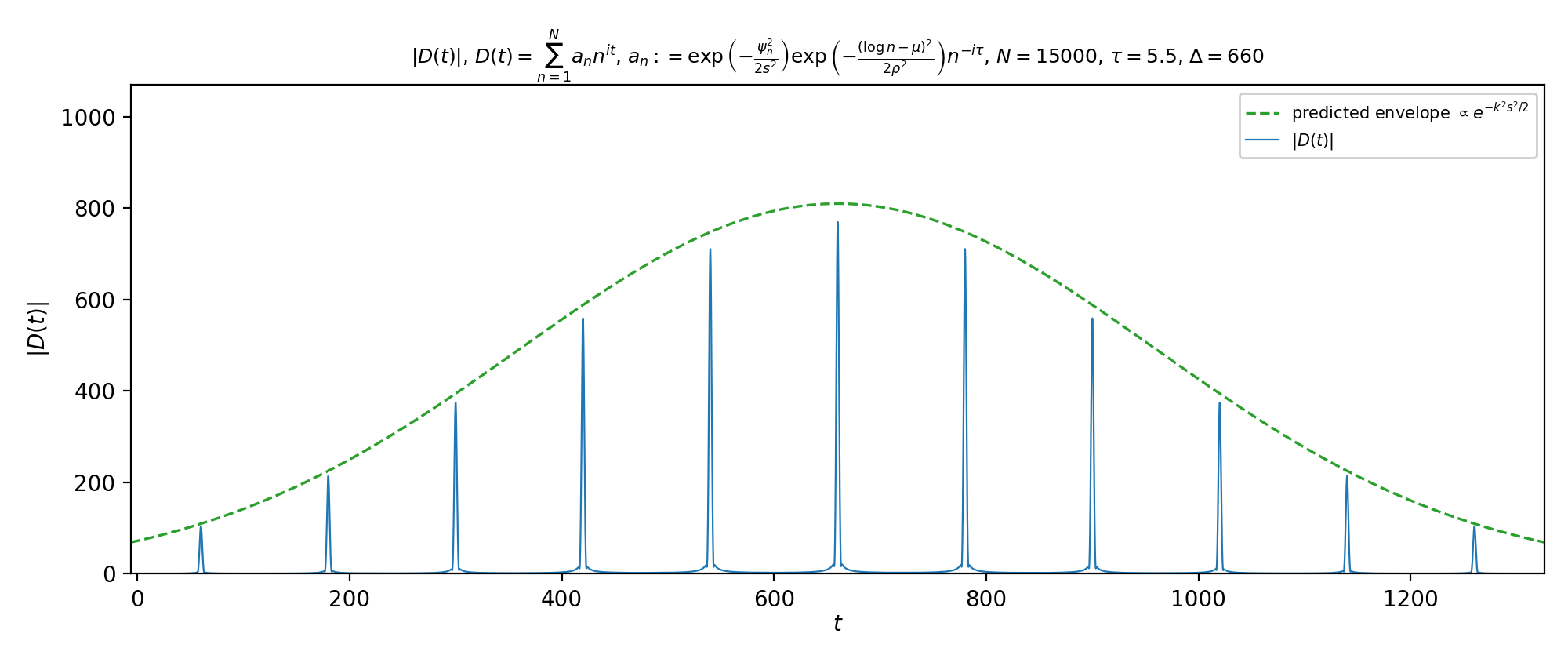}
\caption{An ``arithmetic wave packet'' --- a Dirichlet series with gaussian type coefficients which is concentrated on an arithmetic progression.}\label{arithmetic}
\end{figure}

Unexpectedly, it was shown by Bourgain \cite{bourgain-dirichlet} using Kakeya counterexamples that the original conjecture of Montgomery was too strong to be true!  The idea is ultimately analogous to that of Fefferman's construction, but with the geometric shapes of rectangles or tubes replaced with arithmetic progressions.  Some straightforward calculations show that by choosing the magnitude and phase of coefficients $a_n$ carefully, one could make the Dirichlet polynomial \eqref{dirichlet-polynomial} large when $t$ was close to an arithmetic progression, which is the arithmetic analogue of a ``wave packet'': see Figure \ref{arithmetic}.  By superimposing many such constructions together, one could then create a Dirichlet polynomial whose values became concentrated on a union of such arithmetic progressions, each with different step sizes.  The key observation is then that the type of long thin rectangles or tubes appearing in the Kakeya conjecture can be interpreted as arithmetic progressions if one discretizes space appropriately; see Figure \ref{disc} for an illustration of this idea.

\begin{figure}
\begin{tikzpicture}[scale=1]

  \draw[step=1cm, lightgray, thin] (0,0) grid (10,10);

  \foreach \c/\r in {3/1, 4/3, 5/5, 6/7, 7/9} {
      \pgfmathtruncatemacro{\ybottom}{10-\r}
      \pgfmathtruncatemacro{\ytop}{11-\r}
      \pgfmathtruncatemacro{\xleft}{\c-1}
      \fill[yellow!50] (\xleft, \ybottom) rectangle (\c, \ytop);
  }

  \begin{scope}[shift={(4.5, 5.5)}, rotate=-63.435]
      \draw[thick, blue, fill=blue!20, opacity=0.4] (-5.25, -0.65) rectangle (5.25, 0.65);
  \end{scope}

  \foreach \r in {1,...,10} {
      \foreach \c in {1,...,10} {
          \pgfmathtruncatemacro{\val}{(\r-1)*10 + \c}
          \pgfmathtruncatemacro{\y}{10-\r}
          \node[font=\sffamily] at (\c-0.5, \y+0.5) {\val};
      }
  }

  \draw[thick] (0,0) rectangle (10,10);

\end{tikzpicture}
\caption{By discretizing and enumerating space appropriately, one can convert rectangles (or tubes) into arithmetic progressions.}\label{disc}
\end{figure}

By combining this idea with the Besicovitch--Perron construction, Bourgain was able to show that the strongest forms of the Montgomery conjecture were false.  However, because the compression of that construction was only logarithmic in size, this did not rule out a weaker form of the Montgomery conjecture that would still be strong enough to imply the Lindel\"of hypothesis.  To contradict that conjecture would require instead a failure of the Kakeya conjecture.  Thus, the resolution of the Kakeya conjecture would be a necessary first step towards a resolution of the Lindel\"of hypothesis by this approach!

We have discussed here how various conjectures in analysis, partial differential equations, or analytic number theory imply the Kakeya conjecture, thanks to wave packet type constructions.  It turns out that these implications can be partially reversed: progress on the Kakeya conjecture can be used to make progress on the restriction, Bochner-Riesz, and local smoothing conjectures.  A basic reason for this is the existence of \emph{wave packet decompositions}, which allow for very general oscillating functions to be decomposed into a superposition of wave packets.  This does not fully reduce the analysis of such oscillating functions to Kakeya-type problems, because understanding the exact nature of the superposition - for instance, whether it generates constructive interference, destructive interference, or some combination of the two - is highly delicate.  A powerful technique to handle this that has emerged in recent decades is that of \emph{induction on scales}: use a hypothesis such as the restriction conjecture at smaller spatial scales, combined with the wave packet decomposition, to control the behavior of the function at larger scales.  These ideas have been developed by many authors, particularly Bourgain and Wolff, but it is only in recent years, due to the work of Wang and her coauthors, that they are becoming powerful enough to fully resolve many of the above conjectures.  For instance, the local smoothing conjecture was resolved in two dimensions by Guth, Wang, and Zhang in \cite{gwz} by a very clever combination of wave packet decompositions, induction on scales, and a Kakeya-type estimate.  With the recent resolution of the three dimensional Kakeya conjecture, there is now optimism that many of the other conjectures listed above will also be resolved in the near future; already most of the current ``world records'' on these problems are held by Wang and her coauthors.

The Kakeya conjecture is also a close cousin to other important questions in the general area of geometric measure theory, which has to do with understanding how various geometric operations, such as computing all the distances between elements of a given set, interact with various spatial notions of ``size'', such as the Hausdorff dimension.  These questions include the Falconer distance set conjecture \cite{falconer}, the discretized sum-product problem \cite{bourgain}, the Furstenberg set problem \cite{wolff} (recently solved by Ren and Wang \cite{renwang}), and various questions in projection theory.  We will not go into detail in these other questions and how they relate to the Kakeya conjecture, but refer the reader to \cite{wolff}.  

\section{The road to the solution of the Kakeya conjecture}

It is not to difficult to adapt the two-dimensional Kakeya set constructions to three dimensions, giving a compression factor comparable to $\log(1/\delta)$ as before.  But the argument of C\'ordoba for the two-dimensional problem changes in a fundamental way when one adds the extra dimension: when asking for all the different orientations to be separated from each other by an angle of $\delta$, the number of possible orientations to consider jumps from about $1/\delta$ to about $1/\delta^2$, and so in theory the amount of compression jumps accordingly.  Thus, in principle, the compression factor could be as large as $1/\delta$ (times an additional logarithmic factor).

But this is a worst case analysis, assuming that every tube intersects every other tube.  In two dimensions, this is a reasonable enough scenario --- after all, any two lines in the plane will eventually intersect unless they are parallel --- but is unrealistic in three dimensions, as most pairs of lines will simply be skew to each other, rather than intersecting or be parallel.  

This question was investigated by Bourgain in 1991 \cite{bourgain}, who showed that the upper bound of about $1/\delta$ can be improved 
to about $1/\delta^{2/3}$, and formally proposed Conjecture \ref{kakeya-conj}, which (in our informal phrasing) asserts that the compression factor could not exceed $\delta^{-o(1)}$ in the limit $\delta \to 0$.  The next major progress was by Wolff in 1995 \cite{wolff}, who improved the exponent $2/3$ in Bourgain's bound to $1/2$.  The arguments were complicated, but ultimately rested on a simple geometric observation: if two lines in space intersected each other at some point $P$, then any third line that intersected both of those lines at a different location than $P$ would have to lie in the plane spanned by the first two lines.  Thus, if one was considering lines pointing in different directions, only a small number of them could intersect any intersecting pair of lines.

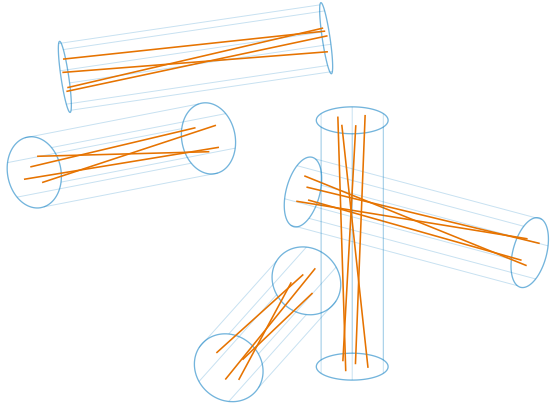
\begin{figure}[h]
\centering
\tdplotsetmaincoords{70}{115} 
\begin{tikzpicture}[tdplot_main_coords, scale=0.95,
   capstyle/.style={draw=cyan!55!blue, line width=0.5pt, opacity=0.55},
   ribstyle/.style={draw=cyan!55!blue, line width=0.3pt, opacity=0.22},
   thinstyle/.style={draw=orange!90!black, line width=0.6pt}]

\def\R{0.50}  
\def\L{1.85}  

\newcommand{\fattube}[3]{%
  \tdplotsetrotatedcoords{#2}{#1}{0}%
  \begin{scope}[shift={#3}]
  \begin{scope}[tdplot_rotated_coords]
    \draw[capstyle] plot[domain=0:360,samples=49,smooth,variable=\t]
        ({\R*cos(\t)},{\R*sin(\t)},{-\L});
    \draw[capstyle] plot[domain=0:360,samples=49,smooth,variable=\t]
        ({\R*cos(\t)},{\R*sin(\t)},{\L});
    \foreach \a in {0,60,...,300}{
      \draw[ribstyle] ({\R*cos(\a)},{\R*sin(\a)},{-\L})
                   -- ({\R*cos(\a)},{\R*sin(\a)},{\L});
    }
    \foreach \ax/\ay/\bx/\by in {0.34/0.18/-0.28/0.40, -0.40/0.26/0.36/-0.10,
                                 0.09/-0.44/0.33/0.29, -0.18/-0.09/-0.40/-0.36}{
      \draw[thinstyle] ({\ax*\R},{\ay*\R},{-\L}) -- ({\bx*\R},{\by*\R},{\L});
    }
  \end{scope}
  \end{scope}
}

\fattube{90}{8}{(0,0,-1.5)}
\fattube{68}{70}{(1.6,-1.5,1.3)}
\fattube{112}{138}{(-1.9,1.4,-0.4)}
\fattube{42}{205}{(1.5,2.0,0.5)}
\fattube{100}{292}{(-1.7,-1.9,1.6)}

\end{tikzpicture}
\caption{Stickiness: the thin tubes (orange) can be organized into a much smaller number of fat tubes (blue), in such a way that thin tubes pointing in nearby directions also lie physically close to each other, and hence are contained in a common fat tube.}\label{sticky}
\end{figure}

\begin{figure}[h]
\centering
\tdplotsetmaincoords{72}{118} 
\begin{tikzpicture}[tdplot_main_coords, scale=0.8,
   planestyle/.style={fill=green!50!teal, fill opacity=0.13,
                      draw=green!45!teal, draw opacity=0.45, line width=0.3pt},
   hairstyle/.style={draw=orange!90!black, line width=0.6pt},
   dotstyle/.style={fill=red!70!black, draw=none}]

\def\Ls{2.6}    
\def\Lh{2.6}    
\def\PA{2.55}   
\def\PB{0.85}   

\foreach \al/\z in {60/-1.7, 105/0.0, 150/1.7}{
  \draw[planestyle]
      ({ \PA*cos(\al)},{ \PA*sin(\al)},{\z-\PB}) --
      ({ \PA*cos(\al)},{ \PA*sin(\al)},{\z+\PB}) --
      ({-\PA*cos(\al)},{-\PA*sin(\al)},{\z+\PB}) --
      ({-\PA*cos(\al)},{-\PA*sin(\al)},{\z-\PB}) -- cycle;
}

\foreach \al/\z in {60/-1.7, 105/0.0, 150/1.7}{
  \foreach \be in {72,81,90,99,108}{
    \draw[hairstyle]
      ({ \Lh*sin(\be)*cos(\al)},{ \Lh*sin(\be)*sin(\al)},{\z+\Lh*cos(\be)}) --
      ({-\Lh*sin(\be)*cos(\al)},{-\Lh*sin(\be)*sin(\al)},{\z-\Lh*cos(\be)});
  }
}

\draw[hairstyle] (0,0,{-\Ls}) -- (0,0,{\Ls});

\foreach \z in {-1.7,0.0,1.7}{
  \draw[dotstyle] (0,0,\z) circle (1.6pt);
}

\end{tikzpicture}
\caption{Planiness, as seen in a ``hairbrush'': a single tube --- the vertical ``stem'' running through all three marked points --- together with all the tubes that meet it.  The tubes passing through any one point of the stem (the red dots) form a ``bush'', and planiness asserts that each such bush is essentially contained in a plane (shaded).  Note that this plane is necessarily one that contains the stem itself, since the stem is one of the tubes through that point; but the plane is free to rotate about the stem as one moves from one point of the stem to the next, as depicted here.}\label{plany}
\end{figure}
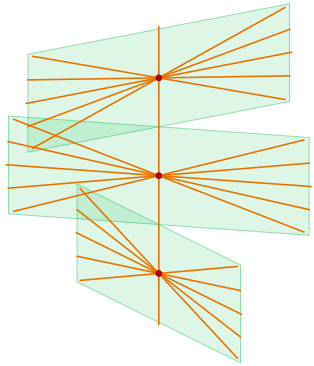

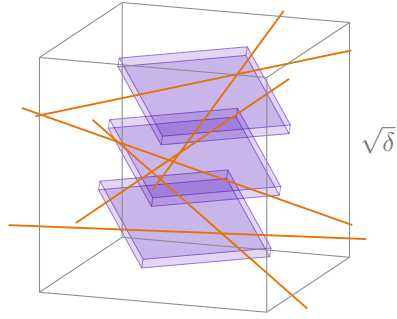
\begin{figure}[h]
\centering
\tdplotsetmaincoords{76}{110} 
\tdplotsetrotatedcoords{40}{22}{0} 
\begin{tikzpicture}[tdplot_main_coords, scale=1.0,
   cubestyle/.style={draw=black!45, line width=0.4pt},
   grainstyle/.style={fill=violet!55!blue, fill opacity=0.16,
                      draw=violet!60!blue, draw opacity=0.55, line width=0.3pt},
   edgestyle/.style={draw=violet!60!blue, draw opacity=0.45, line width=0.3pt},
   tubestyle/.style={draw=orange!90!black, line width=0.7pt}]

\def\h{1.6}    
\def\A{0.9}    
\def\T{0.07}   
\def\Lt{2.4}   

\begin{scope}[tdplot_rotated_coords]
  \foreach \c in {-1,0,1}{
    \foreach \s in {-1,1}{
      \draw[grainstyle] (-\A,-\A,{\c+\s*\T}) -- ( \A,-\A,{\c+\s*\T}) --
                        ( \A, \A,{\c+\s*\T}) -- (-\A, \A,{\c+\s*\T}) -- cycle;
    }
    \foreach \sx/\sy in {-1/-1, 1/-1, 1/1, -1/1}{
      \draw[edgestyle] ({\sx*\A},{\sy*\A},{\c-\T}) -- ({\sx*\A},{\sy*\A},{\c+\T});
    }
  }
\end{scope}

\begin{scope}[tdplot_rotated_coords]
  \foreach \c/\th/\s in {-1/18/0.25, -1/78/-0.30, 0/125/0.20,
                          0/48/-0.25, 1/100/0.30, 1/140/-0.20}{
    \draw[tubestyle] ({-\Lt*cos(\th)-\s*sin(\th)},{-\Lt*sin(\th)+\s*cos(\th)},\c) --
                     ({ \Lt*cos(\th)-\s*sin(\th)},{ \Lt*sin(\th)+\s*cos(\th)},\c);
  }
\end{scope}

\foreach \sa in {-1,1}{\foreach \sb in {-1,1}{
  \draw[cubestyle] (-\h,{\sa*\h},{\sb*\h}) -- (\h,{\sa*\h},{\sb*\h});
  \draw[cubestyle] ({\sa*\h},-\h,{\sb*\h}) -- ({\sa*\h},\h,{\sb*\h});
  \draw[cubestyle] ({\sa*\h},{\sb*\h},-\h) -- ({\sa*\h},{\sb*\h},\h);
}}
\node[black!55, font=\footnotesize, right] at (-\h,\h,0) {$\sqrt{\delta}$};

\end{tikzpicture}
\caption{Graininess: the tubes in an extremal Kakeya configuration tend to fill out ``grains''.  Different formulations of graininess use different dimensions of grains; in this picture, grains of dimensions $\delta \times \sqrt{\delta} \times \sqrt{\delta}$ are depicted.  This property is closely related to planiness.}\label{grainy}
\end{figure}

The exponent of $1/2$ proved extremely resistant to further improvements. In 2000, Katz, {\L}aba, and the author wrote a 64 page paper \cite{klt} in which we were able to lower the exponent $1/2$ incrementally to $1/2 - 10^{-10}$.  This was only a small improvement in the exponent, but the argument revealed the shape of the ``enemy'' blocking the Kakeya conjecture from being resolved: a configuration of lines in space resembling in many key respects an object from quantum mechanics known as the \emph{Heisenberg group}
\begin{equation}\label{heisenberg}
  \{ (z_1,z_2,z_3) \in \C^3 : \mathrm{Im}(z_1 \overline{z_2}) = \mathrm{Im} z_3 \}.
\end{equation}
This is not exactly a Kakeya set for several reasons, the most obvious one being that it lives in $\C^3$ rather than $\R^3$.  Nevertheless, it is Kakeya-like in many ways; in particular, it contains about as many complex lines as one would expect a (complex analogue of a) Kakeya set to contain, while remaining significantly smaller than the entire three-dimensional space $\C^3$.  An analysis of this group suggested that the tubes in a bad Kakeya configuration should obey three key additional properties:
\begin{itemize}
    \item \emph{stickiness}, which roughly speaking asserts that tubes that are oriented in nearby directions must also be physically close to each other, so that the thin tubes can be efficiently organized into ``fat tubes'' as in Figure \ref{sticky};
    \item \emph{planiness}, which roughly speaking asserts that the tubes through any given point were essentially coplanar to each other (thus giving the configuration a sort of ``tangent plane'' at every point), as in Figure \ref{plany}; and
    \item \emph{graininess}, which roughly speaking asserts that the configuration contained medium-sized rectangular slabs (or ``grains''); see Figure \ref{grainy}.
\end{itemize}
With significant effort, the authors of \cite{klt} were able to show that such a configuration was not compatible with the tubes in the configuration pointing in all possible directions when the compression exponent was close to $1/2$.  

Stickiness turns out to be the critical property, allowing for a very useful \emph{self-similar} structure of potential Kakeya counterexamples to surface, in which the portion of a Kakeya configuration inside a ``fat tube'' behaved very much like a rescaled version of a Kakeya configuration at a coarser scale.  This enabled the aforementioned \emph{induction on scales} strategy to become effective, and in particular led to the other two properties of planiness and graininess.  

The author and Katz continued to work on various aspects of the problem for many years, but eventually declared defeat in 2014, laying out a potential road map to solving the problem in \cite{blog} if the properties of stickiness, planiness, and graininess could be justified.  However, there were many obstacles to making this map a rigorous reality, most notably the absence of any mechanism for forcing stickiness to occur.

Katz and Zahl were able to reprove (a slightly improved version of) the results in \cite{klt} in 2017 by a different method \cite{zahl}, and more progress was made in dimensions higher than three; but the three-dimensional Kakeya conjecture remained stubbornly out of reach, until the spectacular complete resolution of the problem by Wang and Zahl in a striking sequence of papers \cite{wang-1, wang-2, wang-3} in 2025 and 2026.

The arguments of Wang and Zahl are described in detail in other places, such as \cite{guth}.  We will just mention a few of the ingredients used here.  Following the road map, the argument splits into two major substeps: a proof of the Kakeya conjecture in the ``sticky'' case, and then a clever reduction of the general case to this sticky case.  Here, the road map was used to obtain the additional properties of planiness and graininess, but Wang and Zahl were also able to extract a very strong ``conjugation'' structure buried inside the Kakeya configuration, behaving very similarly to how complex conjugation $z \mapsto \overline{z}$ is buried inside the Heisenberg group \eqref{heisenberg}.  For instance, in the complex numbers, the expression $z \overline{w} + w \overline{z}$ is necessarily a real number for any $z,w \in \C$; the analogous conjugation operation $s \colon \R \to \R$ arising from a sticky Kakeya configuration also restricts $s(x) y + s(y) x$ to a small set of values for $x,y \in \R$.  Such an operation, as it turns out, can also be used to build counterexamples to another type of statement in geometric measure theory, known as a \emph{projection theorems}.  Wang and Zahl were able to leverage some deep recent progress on such projection theorems or the discretized sum-product conjecture.  However, recent deep progress on such theorems by Orponen, Shmerkin, and Wang \cite{orponen, osw, sw} could then be used to obtain a contradiction; this was a major ingredient that was missing from the original road map.

The reduction to the sticky case is completely original to Wang and Zahl.  At the time of the original road map, the conventional wisdom was that the induction on scales method was only powerful in the sticky case; every attempt to try to control a non-sticky Kakeya configuration by assuming a Kakeya type hypothesis at coarser scales led to losses that rapidly accumulated as one iterated to finer and finer scales.  However, in recent years the power of the induction on scales method has been greatly enhanced, to cover both sticky and non-sticky cases, for instance in the striking resolution by Ren and Wang \cite{renwang} of the Furstenberg set problem.  The means to achieve this are too technical to explain here, save to say that the precise statement one is performing induction on must be carefully tailored to the geometry of the configuration.  The analysis of every non-sticky configuration can then be reduced to the analysis of various rescaled portions of the Kakeya set which, upon iterating enough times, either become sticky, or degenerate to a simpler type of configuration that can be handled by other means, thus establishing the Kakeya conjecture in full generality.

\section{Further directions}

The work of Wang and her co-authors open up a general strategy to handle many other problems in geometric measure theory, harmonic analysis, and beyond, by first isolating the underlying Kakeya-type problem that is controlling the geometry, then treating the sticky case of that problem (as well as cases at the opposite end of the spectrum, such as ``well-spaced'' constructions), and finally to combine these cases with an induction on scales argument to handle the general case.  There is now optimism that many of the other conjectures mentioned in this article, such as the restriction conjecture, the Bochner-Riesz conjecture, and the local smoothing conjecture, will see further progress, or even complete resolution in the near future.  (The Montgomery conjecture, which has less geometry embedded in it, may be particularly challenging, but perhaps even this conjecture will see some further progress.)  It is an exciting time for the field!

\section{Acknowledgments}

AI assistance was used to generate the images in this article, to perform literature search, and to autocomplete text.

\bibliographystyle{amsplain}

\end{document}